\documentstyle[A4,12pt,amssymb,amsfonts]{article}
\begin{document}
\author{Nicusor Dan\footnote{This paper was realized
in the framework of the contract CEx06-10/25.07.06. It was also
supported by the contract CNCSIS 1338/2007.}}
\title{Multitransgression and regulators}
\date{}
\maketitle
\parindent=0pt
\parskip=5pt
\newtheorem{theor}{Theorem}[section]
\newtheorem{prop}[theor]{Proposition}
\newtheorem{cor}[theor]{Corollaire}
\newtheorem{lemma}[theor]{Lemma}
\newtheorem{defin}[theor]{Definition}
\newtheorem{conj}[theor]{Conjecture}
\newtheorem{rem}[theor]{Remark}
\newtheorem{nota}[theor]{Notations}
\def\vs{\vskip 4mm}
\def\gl{G^{\lambda}}
\def\ol{\Omega^\lambda}
\def\dl{{\Delta}^{\lambda}}
\def\gly{G^{\lambda}_Y}
\def\oly{\Omega^{\lambda}_Y}
\def\dly{{\Delta}^{\lambda}_Y}
\def\l{\lambda}
\def\eq{equation}
\def\lu{{\lambda}_1}
\def\ld{{\lambda}_2}
\def\lt{{\lambda}_3}
\def\glu{G^{\lu}}
\def\gld{G^{\ld}}
\def\glt{G^{\lt}}
\def\olu{\Omega^\lu}
\def\old{\Omega^\ld}
\def\olt{\Omega^\lt}
\def\dlu{{\Delta}^{\lu}}
\def\dld{{\Delta}^{\ld}}
 \def\dlt{{\Delta}^{\lt}}
\def\ddl{\frac{\partial}{\partial\l}}
\def\ddlu{\frac{\partial}{\partial\lu}}
\def\ddld{\frac{\partial}{\partial\ld}}
\def\ra{\rightarrow}
\def\ci{{\ C}^\infty}
\def\tix{\tilde{X}}
\def\pii{\pi^{-1}(y)}
\def\curg{courant de Green }
\def\fd{forme diff{\'e}rentielle }
\def\pix{\pi:\tix \ra X}
\def\sl{||s||^{2\l}}
\def\piy{\pi^{-1}(Y)}
\def\espc{espace analytique complexe de dimension finie }
\def\va{\varphi}
\def\ep{\epsilon}
\def\pa{\partial}
\def\al{\alpha}
\def\n{\nabla}
\def\m{\nabla'}
\def\p{\nabla''}
\def\q{\nabla^2}
\def\pc{\prod_k}
\def\mn{\{1,\cdots ,m+n\}}
\def\lz{(\l_l)_l =0}
\def\bw{\bigwedge}
\def\be{\beta}
\def\cd{\cdot}
\def\cds{\cdots}
\def\miu{\mu}
\def\ga{\gamma}
\def\ro{\rho}
\def\xt{X\times T}
\def\de{\delta}
\def\om{\omega}
\def\ul{U_{(\l_l)}}
\def\ri{\rightline{$\Box$}}
\def\rid{\rightline{$\Box$}}
\def\dx{{\cal D}(X)}
\def\dxx{{\cal D}'(X)}
\def\Om{\Omega}
\def\ho{{\cal H}(\Om)}
\def\dih{{\cal D}'_h(X,\Om)}
\def\R{{\Bbb{R} } }
\def\Z{{\Bbb{Z} } }
\def\P{{\Bbb{P} } }
\def\N{{\Bbb{N} } }
\def\Q{{\Bbb{Q} } }
\def\O{{\cal O}}
\def\D{{\cal D}}
\def\A{{\cal A}}
\def\E{{\cal E}}
\def\C{{\Bbb{C}}}
\def\widebar{\overline}

\begin{abstract}
In a parallel way to the work of Wang, we define higher order
characteristic classes associated with the Chern character,
generalizing the work of Bott-Chern and Gillet-Soul\'e on
secondary characteristic classes. Our formalism is simplicial and
the computations are easier. As a consequence, we obtain the
comparison of Borel and Beilinson regulators and an explicit
formula for the real single-valued function associated with the
Grasmannian polylogarithm.

\end{abstract}

\section{Local multitransgression formula}
\label{1000}

Let $X$ be a complex manifold. We denote by
$A(X)=\bigoplus_{p,q}A(X)^{p,q}$ the space of complex differential
forms on $X$. Let $T$ be a real manifold. We denote by $d_X = d' +
d''$ and $d_T$ the differentiation operator on differential forms
on $X \times T$ along $X$ and $T$, respectively.

Let $E$ be a holomorphic vector bundle on $X$. We denote also by
$E$ the inverse image of $E$ on $X \times T$ by the projection $X
\times T \to X$. Let $V=A(\xt, E)$ be the global $\ci$ sections of
$E$ with coefficients differential forms. We denote by $[,] $ the
supercommutator on ${\mathrm End}_\C(V)$ associated to the total
differential degree of differential forms on $X \times T$. We
endow E with a $\ci$ metric $<,>$ on $X \times T$, which is
hermitian restricted to each $X \times \{ t \}, t \in T$. Let
$\nabla = \nabla (t):A^0(X,E) \to A^1(X,E)$ be the holomorphic
unitary connection for this metric, for each point $t \in T$. We
can see $\nabla$ as an application $\nabla :A^0(X\times T, E)\to
A^1(X\times T, E)$ and we can extend it in the usual manner to an
application $\nabla :A^i(X\times T, E)\to A^{i+1}(X\times T, E)$.
We denote by $\nabla',\nabla''$ the holomorphic and the
anti-holomorphic components of $\nabla$ with respect to $X$.

Let $r$ be a positive integer and $\varphi$ a symmetric $\C$
multilinear map on $M_N(\C)^r$, invariant by the adjoint action of
$GL_N(\C)$. We can associate to $\varphi$ a characteristic class
$\alpha := \varphi (\nabla^2, \cdots, \nabla^2)$. If $\varphi =
CH_r$, where $CH_r (A_1,\cdots,A_r)=(-1)^r (r!)^{-1} (r!)^{-1}
\sum_{\sigma\in S_r} Tr(A_{\sigma 1} \cdots A_{\sigma r})$, we
obtain $\alpha = ch_r$, the $r-$th Chern character. We define the
number operator $N\in A^1 (X\times T,{\mathrm End} (E))$ by the
formula $d_T <e,f>=<e,Nf>,$ for any two sections $e,f \in A^0(\xt
, E)$. Locally, if $h$ denotes the hermitian metric expressed in a
base of local sections of $E$, $N=h^{-1}d_T h$. We introduce, for
$1\leq n \leq 2r-1$, the local $n-$transgressed forms of the
characteristic class $\alpha$
$$\alpha^{(n)}=\sum_{\stackrel{k,p,q\geq 0,k+p+q+1 \leq r}{2k+p+q+1=n}} (-1)^{(k+q)}
\frac{(k+p)! (k+q)!) (k+p+q+1)!}{k!p!q!}$$
$$\cdot \varphi (\nabla^{2^{ < r-k-p-q-1 >}}, N, 1/2[N,N]^{< k >}, [\nabla'',N]^{< p >}),
[\nabla',N]^{< q>}),$$ where $X^{< m >}$ means $X,\cdots,X$ ($m$
times).

\begin{theor}
\label{theor1}
$$d'_X d''_X \alpha^{(1)}+d_T \alpha =0,$$
\begin{\eq}
d_X \alpha^{(n)}+ n d_T \alpha^{(n-1)} = \varphi (\nabla^{2^{
< r-n
>}},[\nabla'',N]^{< n >} - (-1)^n [\nabla',N]^{< n >})\ \  (n\geq 2).
\label{100}
\end{\eq}
\end{theor}

\section{The definition of the regulator}
\label{2000}

Let $X$ be a smooth projective complex variety. We denote by
$A_{\R}^{\cdot}(X) $ the set of complex differential forms on $X$
invariant by complex conjugation. We define
$A_{\R}^{\cdot}(X)(r):=(2i \pi)^r A_{\R}^{\cdot}(X)\subset
A^{\cdot}(X)$. We recall from [W], [B1] the definition of the
complex $\D^{\cdot}(X,r)$ computing the real Deligne cohomology:
$${\D}^{n}(X,r)=A_{\R}^{n-1}(X)(r-1)\bigcap \bigoplus_{p+q=n-1,
p<r,q<r}A^{p,q}(X)\ \ \  {\mathrm if}\  n\leq 2r-1,$$
$${\D}^{n}(X,r)=A_{\R}^{n}(X)(r)\bigcap \bigoplus_{p+q=n,
p\geq r,q\geq r}A^{p,q}(X) \ \ \ \ \ \ {\mathrm if}\  n\geq 2r,$$
with the differential of the element $x\in {\D}^{n}(X,r)$ given by
$$d_{\D}x=dx \ \ \ \ \ \ {\mathrm if}\  n\geq 2r,$$
$$d_{\D}x=-2d'd''x \ \ \ \ \ \ {\mathrm if}\  n=2r-1,$$
$$d_{\D}x=-\pi (dx)\ \ \ \ \  {\mathrm if}\  n<2r-1,$$
where $\pi (d)$ is the restriction of $d$ to the $p\leq n, q\leq
n$ pieces of the bigrading.

We define the simplicial set $T_{\cdot}(X)$ of isomorphism classes
of hermitian vector bundles on $X$. An element $\bar{E}
=(E_0,...E_n;h_0,...,h_n)\in T_n(X)$ consist of vector bundles
$E_0,...E_n$ on $X$ endowed with hermitian metrics $h_0,...,h_n,$
and isomorphisms $\varphi_i : E_0 \to E_i$ for each $i=1,...,n$.

We denote by ${\Delta}^n$ the standard $n-$simplex ${\Delta}^n=\{
(t_0,...t_n), t_i\geq 0, \sum_{i=0}^n t_i =1\}$. We associate to
an element $\bar{E}\in T_n(X)$ the hermitian vector bundle
$\bar{\E} =(\E, h_t)$ on $X \times {\Delta}^n$, where the vector
bundle $\E$ is the inverse image of the vector bundle $E_0$ by the
projection $X \times {\Delta}^n \to X,$ and the metric $h_t$ on
 $\E$ restricted to $X\times \{ (t_0,...,t_n)\}$ is $t_0
h_0+t_1 \varphi_1^{*}h_1+...t_n (\varphi_n ... \varphi_1)^{*}h_n$.
We define:
$${\mathrm ch}_r^{(n)}(\bar{E}) =\frac{1}{2n!} \int_{\Delta^n} ch_r^{(n)}(\bar{\E}).$$

We observe that ${\mathrm ch}_r^{(n)}(\bar{E}) \in
{\D}^{2r-n}(X,r)\subset A^n(X)$. We extend by linearity the
application ${\mathrm ch}_r^{(n)}$ to an application ${\mathrm
ch}_r^{(n)}: \Z T_n (X) \to  {\D}^{2r-n}(X,r)$. Theorem
$\ref{theor1}$ and Stokes formula imply $d_{\D} {\mathrm
ch}_r^{(n)}(\E)={\mathrm ch}_r^{(n-1)}(\pa \E)$ (the right term of
equation $(\ref{100})$ disappears because of its complex
bigrading), i.e. ${\mathrm ch}_r^{(\cdot)}: \Z T_{\cdot} (X) \to
{\D}^{2r-\cdot}(X,r)$ is a morphism of complexes.

If we fix a metric $h$ on the trivial rank $N$ vector bundle $1^N$
on $X$ and if we consider only those elements of the form $\bar{E}
= (1^N,...,1^N;h,...,h)$ in $T_n(X)$, we obtain an inclusion of
simplicial sets $B.GL_N(X)\to T_{\cdot}(X)$. We obtain therefore a
morphism:
\begin{\eq}
\label{200}
 {\mathrm ch}_r^{(\cdot)}: \Z B.GL_N(X)\to
{\D}^{2r-\cdot}(X,r).
\end{\eq}

\begin{prop}
\label{prop1} The morphism $(\ref{200})$ is, up to homotopy,
independent of the choice of the metric $h$ and compatible with
the inclusion $B.GL_N(X)\subset B.GL_{N+1}(X).$
\end{prop}

So we obtain a morphism $H_n(B.GL_{\infty}(X)) \to
H^{2r-n}_{\D}(X,\R(r))$. If $X={\mathrm Spec} \C$, we compose this
morphism with the classical morphism $K_n(\C) =
\pi_n(B.GL_{\infty}(\C)^+) \to
H_n(B.GL_{\infty}(\C)^+)=H_n(B.GL_{\infty}(\C))$ and we obtain a
morphism
\begin{\eq}
\label{300}
 {\mathrm ch}_r^{(n)} : K_n(\C) \to H^{2r-n}_{\D}(\C,\R(r)).
\end{\eq}

\begin{theor}
\label{theor2} The morphism $(\ref{300})$ coincides with the
Beilinson regulator.
\end{theor}
\par{\bf Proof: }Let $S_{\cdot}(X)$ be the simplicial set computing
the Waldhausen K-theory of $X$. Wang ([W]) constructs an explicit
morphism of complexes $\Z S_{\cdot}(X) \to {\D}^{2r+1-\cdot}(X,r)$
and Burgos-Wang ([BW]) proved that the composition
$K_n(X):=\pi_{n+1}(S_{\cdot}(X))\to H_{n+1}(\Z S_{\cdot}(X))\to
H^{2r-n}_{\D}(X,\R(r))$ is the Beilinson regulator. There is a
canonical map $\Sigma B_{\cdot}GL_N(X)\to S_{\cdot}(X)$, so we
deduce a morphism ${\mathrm ch}_r^{(\cdot)}: \Z B.GL_N(X)\to
{\D}^{2r-\cdot}(X,r)$. We prove that this morphism coincides with
the morphism $(\ref{200})$ by using an explicit homotopy on the
bisimplicial complex $\Sigma B_{\cdot \cdot}GL_N(X)$.

\section{Comparison of Borel and Beilinson regulators}
The morphism $(\ref{200})$ applied to an element of
$B_{2r-1}GL_N(X)$ written in the homogeneous form
$(g_0,\cdots,g_{2r-1})$ gives
\begin{\eq}
h_t=\sum_{i=0}^{2r-1} t_i g_i h \bar{g_i}^t,
 \label{400}
\end{\eq}
\begin{\eq}
\label{500} {\mathrm
ch}_r^{(2r-1)}(g_0,\cdots,g_{2r-1})=-\frac{(r-1)!}{2(2r-1)!}\int_{\Delta^n}
Tr(h_t^{-1} d_T h_t)^{2r-1}.
\end{\eq}
So the morphism $H_{2r-1}(B.GL_N(\C)) \to H^{1}_{\D}(\C,\R(r))=\R
(r-1)$ giving rise to the Beilinson regulator can be explicited by
these formulas. On the other side, the Borel regulator comes ([H])
from the explicit morphism $H_{2r-1}(B.GL_N(\C)) \to \R $ given by
$b'_r(g_0,\cdots,g_{2r-1})={\mathrm cst}\int_{\Delta^n}
Tr(h_t^{-1} d_T h_t)^{2r-1}$ where $h_t=\sum_{i=0}^{2r-1} t_i g_i
\bar{g_i}^t$ ($h=1$ in $(\ref{400})$). Taking care of the
normalizations, we obtain the theorem of [B2]
\begin{theor}
The normalized Borel regulator is twice the Beilinson regulator.
\end{theor}

\section{Explicit presentation of the real period attached to the
grasmmanian polylogarithm}

The integral $(\ref{500})$ still converges if the metric $h$ is
degenerate of rank $N-r+1$ and we prove that we obtain the same
morphism $H_{2r-1}(B.GL_N(X)) \to H^{1}_{\D}(X,\R(r))$. When $N=r$
we can use a metric of rank one $h=v\bar{v}^t$ for a everywhere
nonvanishing section $v\in H^0 (X,1^N)$. Let $v_i=g_i v$. We
denote by ${{\cal I}}^s$ the set of subsets of $\{
0,\cdots,2r-1\}$ of cardinality $s$ and for $I=\{i_1,\cdots,i_s\}
\in {\cal I}^s$ we denote $t_I=\prod_{i\in I} t_i$ and
$v_I=v_{i_1},\cdots,v_{i_s}$. The morphism $(\ref{200})$ is

\begin{theor}
$${\mathrm
ch}_r^{(2r-1)}(g_0,\cdots,g_{2r-1})=-\frac{(r-1)!}{2(2r-1)!}\cdot$$
\begin{\eq}
\cdot \int_{\Delta^n}\frac{\sum_{\stackrel{I_1,\cdots,I_{2r-1}\in
{\cal I}^{r-1}}{0\leq i_1,\cdots,i_{2r-1}\leq
2r-1}}\prod_{j=1}^{2r-1} t_{I_j}dt_{i_j}
\widebar{det(v_{I_j},v_{i_j})}det(v_{I_j},v_{i_{j+1}})}{(\sum_{I\in
{\cal I}^r} t_I |det(v_I)|^2)^{2r-1}}. \label{600}
\end{\eq}
\end{theor}
Goncharov ([G]) associates to $2r$ non-zero vectors
$v_0,\cdots,v_{2r-1}$ in ${\C}^r$ a $\Z$ mixed Hodge structure
${{\cal G}}(v_0,\cdots,v_{2r-1})$ of Hodge-Tate type, the maximal
period of which is the Grasmannian polylogarithm.
\begin{conj}
The maximal period of the $\R$-MHS attached to ${\cal
G}(v_0,\cdots,v_{2r-1})$ is the right-hand side of $(\ref{600})$.
\end{conj}
\begin{theor}
The conjecture is true for $r\leq 3$.
\end{theor}
For $r=2$ we have some simplifications in the right-hand side of
$(\ref{600})$, beginning with the nice classical type
\begin{prop}
The following function is antisymmetric with respect to
permutations in the variables $v_i$:
$$f(v_0,v_1,v_2,v_3)=Im(det(v_0,v_1)det(v_2,v_3)\widebar{det(v_0,v_3)}\widebar{det(v_1,v_2)}).$$
\end{prop}
We obtain a new presentation of the Bloch-Wigner dilogarithm $\R
Li_2$:
\begin{theor}
Let $v_0,v_1,v_2,v_3$ be non-zero vectors in ${\C}^2$ and
$r(v_0,v_1,v_2,v_3)$ the cross-ratio of the four points they
represent in ${\P}^1_{\C}$. Then
$$\R
Li_2(r(v_0,v_1,v_2,v_3))=12 i f(v_0,v_1,v_2,v_3)\cdot
\int_{\Delta^3} \frac{dt_1 dt_2 dt_3}{(\sum_{i\neq j}t_i t_j
|det(v_i,v_j)|^2)^2}.$$
\end{theor}

\section{A speculation concerning the Beilinson-Soul\'e conjecture}
It is remarkable that the vanishing of ${\mathrm ch}_r^{(n)}$ for
$n\geq 2r$ is a consequence of the local formalism of the section
$\ref{1000}$ and not of the vanishing of $\D^{n}(X,r)$ for $n\geq
2r$. If instead of the complex $\D^{\cdot}(X,r)$ we had a complex
${\cal B}^{\cdot}(X,r)$ and a theory ${\mathrm ch}_r^{(n)} :
K_n(X) \to H^{2r-n}({\cal B}^{\cdot}(X,r))$ which: a) could be
realized by a formalism as in section $\ref{1000}$ and b) is
injective on $K_n^{[r]}(X)\otimes \Q$ (the weight $r$ piece of
K-theory for the Adams operations), we would have the
Beilinson-Soul\'e conjecture.

\vs \par{\bf References:}

[BC]: R. Bott, S.S. Chern: Hermitian vector bundles and the
equidistribution of zeroes of their holomorphic sections, Acta
Math. 114(1968), p. 71-112

[B1]: J.I. Burgos: Green forms and their product, Duke. Math. J.
75(1994), p. 529-574

[B2]: J.I. Burgos: The Regulators of Beilinson and Borel, CRM
Monogr. Ser., vol. 15, AMS, Providence, RI (2001)

[BW]: J.I. Burgos, S. Wang: Higher Bott-Chern forms and
Beilinson's regulator, Invent. Math. 132(1998), p. 261-305

[G]: A. B. Goncharov: Chow polylogarithms and regulators, Math.
Res. Letters 2(1995), p. 99-114

[GS]: H. Gillet, C. Soul\'e: Characteristic classes for algebraic
vector bundles with hermitian metric, Ann. of Math. 131(1990), p.
163-238

[H]: N. Hamida: Description explicite du regulateur de Borel, C.
R. Acad. Sci. Paris Sr. I Math. 330(2000), p. 169-172

[W]: S. Wang: Higher-order characteristic classes in arithmetic
geometry, Thesis Harvard, 1992

\vs Institute of Mathematics of the Romanian Academy, Calea
Grivitei 21, 010702 Bucharest, Romania

E-mail: ndan@dnt.ro

\end{document}